\documentclass[reqno,centertags, 12pt]{amsart}
\usepackage{amscd}
\usepackage{latexsym}
\usepackage{amsmath}
\usepackage{amsthm}
\usepackage{amsfonts}
\usepackage{amssymb}
\numberwithin{equation}{section}

\newcommand{\vp}{\varphi}
\newcommand{\T}{\partial\mathbb{D}}

\newcommand{\ds}{\displaystyle}
\newcommand{\ol}{\overline}

\newcommand{\bigonsq}{O\left( \frac{1}{n^2} \right)}

\newcommand{\be}{\begin{equation}}
\newcommand{\ee}{\end{equation}}
\newcommand{\ba}{\begin{array}}
\newcommand{\ea}{\end{array}}
\newcommand{\C}{\mathbb{C}}
\newcommand{\D}{\mathbb{D}}

\newtheorem{theorem}{Theorem}[section]

\begin{document}

\title[Generalized Bounded Variation and Inserting Point Masses]
{Generalized Bounded Variation and Inserting Point Masses}
\author[M.-W. L. Wong]{Manwah Lilian Wong}
\thanks{$^*$ MC 253-37, California Institute of Technology, Pasadena, CA 91125.
E-mail: wongmw@caltech.edu }
\date{February 10th, 2007}

\begin{abstract} Let $d\mu$ be a probability measure on the unit circle and $d\nu$ be the measure formed by adding a pure point to $d\mu$. We give a formula for the Verblunsky coefficients of $d\nu$ following the method of Simon.

Then we consider $d\mu_0$, a probability measure on the unit circle with $\ell^2$ Verblunsky coefficients $(\alpha_n (d\mu_0))_{n=0}^{\infty}$ of bounded variation. We insert $m$ pure points $z_j$ to $d\mu_0$, rescale, and form the probability measure $d\mu_m$. We use the formula above to prove that the Verblunsky coefficients of $d\mu_m$ are in the form $\alpha_n(d\mu_0) + \sum_{j=1}^m \frac{\ol{z_j}^{n} c_j}{n} + E_n$, where the $c_j$'s are constants of norm $1$ independent of the weights of the pure points and independent of $n$; the error term $E_n$ is in the order of $o(1/n)$. Furthermore, we prove that $d\mu_m$ is of $(m+1)$-generalized bounded variation - a notion that we shall introduce in the paper. Then we use this fact to prove that $\lim_{n \to \infty} \vp_n^*(z, d\mu_m)$ is continuous and is equal to $D(z, d\mu_m)^{-1}$ away from the pure points.

\end{abstract}
\maketitle

\section{Introduction}

Suppose we have a probability measure $d\mu$ on the unit circle $\T=\{z \in \C: |z|=1\}$. We define an inner product and a norm on $L^2(\T, d\mu)$ respectively by
\begin{eqnarray}
\left\langle f ,g \right\rangle & = & \ds \int_{\T} \ol{f(e^{i \theta})} g(e^{i \theta}) d\mu(\theta) \\
\|f\|_{d\mu} & = & \left( \ds \int_{\T} |f(e^{i \theta})|^2 d\mu(\theta) \right)^{1/2}
\end{eqnarray}

Then we orthogonalize $1, z, z^2, \dots$ to obtain the family of monic orthogonal polynomials associated with the measure $d\mu$, namely, $(\Phi_n(z, d\mu))_{n=0}^{\infty}$. We denote the normalized family as $(\vp_n(z, d\mu))_{n=0}^{\infty}$.

The family of orthogonal polynomials on the unit circle obey the Szeg\H o recursion relation: let $\Phi_n^*(z)=z^n \ol{\Phi_n(1/\ol{z})}$ and $\vp_n^*(z)=\Phi_n^*(z)/\|\Phi_n\|$ (they are often known as the reversed polynomials). Since $\Phi_n(z)$ is the unique $n^{th}$ degree monic polynomial that is orthogonal to $1, z, \dots, z^{n-1}$, $\Phi_n^*(z)$ is the unique polynomial of degree $\leq n$ (up to multiplication by a constant) that is orthogonal to $\{z, z^2. \dots, z^n\}$. Then we note that $\Phi_{n+1}(z) - z \Phi_n (z)$ is a polynomial of degree at most $n$ which is orthogonal to $z, z^2, \dots, z^n$, hence, there exists a constant $\alpha_n$ such that the following holds
\be
\Phi_{n+1}(z) = z \Phi_n (z)- \ol{\alpha_n} \Phi_n^*(z)
\label{eq01}
\ee
$\alpha_n$ is called the $n$-th Verblunsky coefficient. From (\ref{eq01}), we could deduce the following recurrence relation for $\Phi_n^*$
\be
\Phi_{n+1}^*(z) = \Phi_n^*(z) - \alpha_n z \Phi_n(z)
\label{eq02}
\ee

Now we consider the norms of the left hand side and the right hand side of (\ref{eq01}) respectively. First, observe that $\|z\Phi_n\|$ is just $\|\Phi_n\|$. Then note that $\Phi_n^*(z)$ is of degree strictly less than $n+1$, so it is orthogonal to $\Phi_{n+1}$. Besides, $\|\Phi_n\|=\|\Phi_n^*\|$. As a result, we have
\be
\|\Phi_{n+1}\|^2 = (1 - |\alpha_n|^2) \| \Phi_n \|^2 = \ds \prod_{j=0}^{n} (1-|\alpha_j|^2)
\label{normphi}
\ee
This also proves that $\alpha_n \in \D=\{z \in \C: |z|<1 \}$.
From (\ref{eq01}), (\ref{eq02}) and (\ref{normphi}) above, we can deduce the Szeg\H o recursion relations for the normalized families as well
\begin{align}
\vp_{n+1}(z) & = (1-|\alpha_n|^2)^{-1/2}(z \vp_n(z)-\ol{\alpha_n} \vp_n^*(z)) \label{normrec} \\
\vp_{n+1}^*(z) & = (1-|\alpha_n|^2)^{-1/2}(\vp_n^*(z) - \alpha_n z \vp_n(z))
\end{align}

From the arguments above, we see that each non-trivial probability measure on the unit circle $d\mu$ corresponds to a sequence $(\alpha_n(d\mu))_{n=0}^{\infty}$ in $\D^{\infty}$ called the Verblunsky coefficients. In fact, the reverse is also true by Verblunsky's theorem, i.e., any sequence of complex numbers $(a_n)_{n=0}^{\infty} \in \D^{\infty}$ is the family of Verblunsky coefficients of a unique probability measure on the unit circle. Hence, there is a bijective correspondence between $(\alpha_n(d\mu))_{n=0}^{\infty}$ and $d\mu$.

The family of Verblunsky coefficients often gives important information about the measure and the family of orthogonal polynomials, for example, from (\ref{normphi}) we know that $\sum_{j=0}^{\infty} |\alpha_j|^2 < \infty$ implies that
\be
\lim_{n \to \infty} \|\Phi_n\| = \ds \prod_{j=0}^{\infty}(1-|\alpha_j|^2)^{1/2} > 0
\ee
This is a fact that we shall use later in the paper. For a more comprehensive introduction to the theory of orthogonal polynomials on the unit circle, the reader should refer to \cite{simon1, simon2}, or the classic references \cite{geronimus, szego}.
\medskip

\section{Results}

In this paper we are going to prove three results, the first one being the following formula
\begin{theorem} Suppose $d\mu$ is a probability measure on the unit circle and $0<\gamma<1$. Let $d\nu$ be the probability measure formed by adding a point mass $\zeta = e^{i \omega} \in \T$ to $d\mu$ in the following manner
\be
d\nu = (1-\gamma) d\mu + \gamma \delta_{\omega}
\label{dnu}
\ee
Then the Verblunsky coefficients of $d\nu$ are given by
\be
\alpha_n(d\nu) = \alpha_n + \ds \frac{(1-|\alpha_n|^2)^{1/2}}{(1-\gamma) \gamma^{-1}+ K_n(\zeta)} \ol{\vp_{n+1}(\zeta)} \vp_n^*(\zeta)
\label{myformula}
\ee where
\be
K_n(\zeta) = \ds \sum_{j=0}^{n} |\vp_j(\zeta)|^2
\ee and all objects without the label $(d\nu)$ are associated with the measure $d\mu$.
\label{theorem0}
\end{theorem}

Before we state the second result, we need to introduce the notion of \emph{p-generalized bounded variation}, $W_{p}(\zeta_1, \zeta_2, \dots, \zeta_p)$, which is the class of sequences defined as follows \smallskip

\textbf{Definition} We say that a sequence $(\alpha_n)_{n=0}^{\infty}$ is of \emph{p-generalized bounded variation} if each $\alpha_n$ can be decomposed into $p$ components
\be
\alpha_n = \ds \sum_{k=1}^{p} \beta_{n,k}
\label{bvcondition1}
\ee with $\beta_{n,k} \in \mathbb{C}$ and there exist $\zeta_1, \zeta_2, \dots, \zeta_p \in \T$ such that for each $1 \leq k \leq p$
\be
\ds \sum_{n=0}^{\infty} |\zeta_k \beta_{n+1,k} - \beta_{n,k}| < \infty
\label{bvcondition2}
\ee We denote by $W_{p}(\zeta_1, \zeta_2, \dots, \zeta_p)$ the class of sequences $(\alpha_n)_{n=0}^{\infty}$ that satisfy (\ref{bvcondition1}) and (\ref{bvcondition2}).

In particular, when $p=1$ and $\zeta_1=1$, then it becomes the conventional bounded variation. This is why we gave the name $p$-generalized bounded variation.

For the sake of simplicity, we shall write $d\mu \in W_{p}(\zeta_1, \zeta_2, \dots, \zeta_p)$ if the family of Verblunsky coefficients of $d\mu$ is in the class $W_p(\zeta_1, \zeta_2, \dots, \zeta_p)$.

\smallskip

The Szeg\H{o} function, which will be involved in Theorem \ref{mainthm2}, is defined as follows

\textbf{Definition} If $d\mu = w(\theta)\frac{d\theta}{2\pi} + d\mu_s$ and $\sum_{j=0}^{\infty} |\alpha_j|^2 < \infty$,  the Szeg\H{o} function is defined as
\be
D(z) = \exp \left(\ds \frac{1}{4\pi} \ds \int \frac{e^{i \theta} + z}{e^{i \theta} -z} \log w(\theta) d\theta \right)
\label{szegofunction}
\ee The well-known Szeg\"o's Theorem asserts the following equality
\be
\ds \prod_{j=0}^{\infty} (1 - |\alpha_j|^2) = exp\left( \int_{0}^{2\pi} \log(w(\theta)) \frac{d\theta}{2 \pi}\right)
\ee Hence, if $(\alpha_n)$ is $\ell^2$, $\log w(\theta)$ is integrable and $D(z)$ defines an analytic function on $\mathbb{D}$. For a thorough discussion of the Szeg\"o function, the reader may refer to Chapter 2 of \cite{simon1}.
\medskip

Now we are ready to state the other two results in this paper:

\begin{theorem} Let $\zeta_j = e^{i \omega_j} \in \T$, $1 \leq j \leq p$ be distinct. Suppose we have a measure $d\mu$ with $d\mu \in W_p (\zeta_1, \zeta_2, \dots, \zeta_p)$ such that for each $j$, $(\beta_{n,j})_{n=0}^{\infty} \in \ell^2$. The following two results hold \\

(1) For any compact subset $K$ of $\T \backslash\{ \zeta_1, \zeta_2, \dots, \zeta_p\}$,
\begin{equation}
\ds \sup_{n; z \in K} |\Phi_n^*(z)| < \infty
\label{mainthm21}
\end{equation}

(2) The following limits are continuous at $z \not = \zeta_1, \zeta_2, \dots, \zeta_p$
\begin{eqnarray}
\Phi_{\infty}^*(z) & = \ds \lim_{n \to \infty} \Phi_n^*(z) & = D(0)D(z)^{-1} \\
\vp_{\infty}^*(z) & = \ds \lim_{n \to \infty} \vp_n^*(z) & = D(z)^{-1}
\label{mainthm22}
\end{eqnarray} and the convergence is uniform on any compact subset $K  \subset \T \backslash \{\zeta_1, \zeta_2, \dots, \zeta_p\}$. Moreoever, $d\mu_s$ is a pure point measure supported on a subset of $\{\zeta_1, \zeta_2, \dots, \zeta_p \}$.
\label{mainthm2}
\end{theorem}
\medskip

\begin{theorem}\label{mainthm1}
Suppose $d\mu_0 \in W_1(1)$ and $(\alpha_n(d\mu_0))_{n=0}^{\infty} \in \ell^2$. We add $m$ distinct pure points $z_j = e^{i \omega_j}$, $\omega_j \not = 0$, to $d\mu_0$ with weights $\gamma_j$ to form the probability measure $d\mu_m$ as follows
\begin{equation}
d\mu_m = \left(1-\ds \sum_{j=1}^{m} \gamma_j \right) d\mu_0 + \ds \sum_{j=1}^{m} \gamma_j \delta_{\omega_j}
\label{dmumdef}
\end{equation} under the conditions that $0  < \gamma_j$ and $\sum_{j=1}^m \gamma_j < 1$. Then $d\mu_m\in W_{m+1}(1, z_1, z_2, \dots, z_m)$ and
\begin{equation}
\alpha_n(d\mu_m) = \alpha_n(d\mu_0) + \ds \sum_{j=1}^{m} \frac{\ol{z_j}^{n} c_j}{n} + E_n
\label{equation1}
\end{equation} where $c_j  = \ol{z_j} |D(z_j, d\mu_0)|^{2} D(z_j, d\mu_0)^{-2}$ are constants independent of the weights $\gamma_1, \gamma_2, \dots, \gamma_m$ and of $n$; and
\be
E_n = E_n(z_1, z_2, \dots, z_m, \gamma_1, \gamma_2, \dots, \gamma_m)= o\left( \ds \frac{1}{n}\right)
\label{enbound}
\ee  Furthermore, for $z \in \T \backslash \{1, z_1, z_2, \dots, z_m\}$, $\vp_\infty^*(z, d\mu_m)$ is continuous and is equal to $(1-\sum_{j=1}^{m} \gamma_j)^{-1/2} D(z, d\mu_0)^{-1}$.\smallskip
\end{theorem}

\textbf{Remark:} Note that ${d\mu_m}_{a.c.}$ is just $(1-\sum_{j=1}^{m} \gamma_j) {d\mu_0}_{a.c.}$ and that $\int \frac{e^{i \theta}+z}{e^{i \theta} -z} \frac{d\theta}{2\pi} = 1$. Hence, $D(z, d\mu_m)=(1-\sum_{j=1}^{m} \gamma_j)^{1/2} D(z, d\mu_0)$. 

\bigskip

Theorem \ref{mainthm2} is a generalization of the following result of Nevai \cite{nevai1} and Nikishin \cite{nikishin}  which reads
\begin{theorem} Suppose $\sum_{j=0}^{\infty}|\alpha_j|^2<\infty$ and
\begin{equation}
\ds \sum_{j=0}^{\infty}|\alpha_{j+1}-\alpha_j|<\infty
\label{boundedvariation}
\end{equation}
Then, for any $\delta>0$,
\begin{equation}
\sup_{n; \delta<\arg(z)<2\pi-\delta} |\Phi_n^*(z)|<\infty
\end{equation} and away from $z=1$, we have that $\lim_{n\to \infty} \Phi_n^*(z)$ exists, is continuous and equal to $D(0) D(z)^{-1}$. Furthermore, $d\mu_s=0$ or else a pure point at $z=1$.
\label{bookthm}
\end{theorem} The reader may refer to Theorem 10.12.5 of \cite{simon2} for the proof. \\
\medskip

According to Simon \cite{simon2}, the history of the problem is as follows. The earliest work related to adding point masses was done by Wigner-von Neumann \cite{wigner}, where they constructed a potential with an embedded eigenvalue. Later, Gel'fand-Levitan \cite{gelfand} constructed a potential $V$ so that $-\frac{d^2}{dx^2}+V$ has a spectral measure with a pure point mass at a positive energy and was otherwise equal to the free measure. A more systematic approach to adding point masses to a potential was then taken by Jost-Kohn \cite{jost1, jost2}.

Unaware of the Jost-Kohn work and of each other, formulae for adding point masses for orthogonal polynomials on the real line case were found by Uvarov \cite{uvarov} and Nevai \cite{nevai}. They found the perturbed polynomials, and Nevai computed the perturbed recursion coefficients.

Jost-Kohn theory for orthogonal polynomials on the unit circle appears previously in Cachafeiro-Marcell\'an \cite{cm1, cm2, cm3}, Marcell\'{a}n-Maroni \cite{mm}, and Peherstorfer-Steinbauer \cite{ps}. In particular, if $d\nu$ and $d\mu$ are as defined in (\ref{dnu}) above, Peherstorfer-Steinbauer \cite{ps} proved that boundedness of the first and second kind orthonormal polynomials of $d\mu$ at the pure point $\zeta$ implies that $\lim_{n \to \infty} \alpha_n(d\nu) - \alpha_n(d\mu) = 0$, but they did not establish any rate of convergence.

When I proved (\ref{myformula}), I was unaware of the following formula found by Geronimus \cite{geronimus}
\be
\Phi_{n}(z, d\nu) = \Phi_n(z) - \ds \frac{\Phi_n(\zeta) K_{n-1}(z, \zeta)}{(1-\gamma) \gamma^{-1} + K_{n-1}(\zeta, \zeta)}
\label{geronimus}
\ee

Years after Geronimus proved (\ref{geronimus}), a similar formula for the real case was rediscovered by Nevai \cite{nevai}, and the same formula for the unit circle case was rediscovered by Cachafeiro-Marcellan \cite{cm3}. Unaware of Geronimus' result and the fact that Nevai's result also applies to the unit circle, Simon reconsidered this problem and proved formula (\ref{mainformula}) independently using a totally different method (see Theorem 10.13.7 in \cite{simon2}). However, a more useful form of his result (see formula (\ref{mainformula}) in Section \ref{proof0}) is disguised in his proof and it lays the foundation to Theorem \ref{theorem0}.\smallskip

In addition to Nevai, Uvarov and Simon's result mentioned above, we use Pr\"ufer variables as the main tool to prove that $\lim_{n \to \infty} \Phi_n^*(z)$ exists. Pr\"ufer variables are named after Pr\"ufer \cite{pruefer}. Their initial introduction in the spectral theory of orthogonal polynomials on the unit circle was made by Nikishin \cite{nikishin} with a significant follow up by Nevai \cite{nevai1}. Both \cite{nevai1} and \cite{nikishin} had results related to Theorem \ref{bookthm} and they arrived at the result by essentially the same proof. Later, Pr\"ufer variables were used as a serious tool in spectral theory by Kiselev-Last-Simon \cite{kiselev} and Last-Simon \cite{last}.
\medskip

Most recently, in \cite{simon1} (Example 1.6.3, p. 72) Simon considered the measure $d\nu$ with one pure point
\be
d\nu = (1-\gamma) \frac{d\theta}{2\pi} + \gamma \delta_{0}
\label{dnudef}
\ee
He proved that the $n$-th degree orthogonal polynomial of $d\nu$ is as follows
\be
\Phi_n(z)=z^n - \ds \frac{\gamma}{1+(n-1) \gamma}(z^{n-1}+z^{n-2} + \dots + 1)
\ee and since $\alpha_n = -\ol{\Phi_{n+1}(0)}$, 
\be
\alpha_n(d\nu) = \ds \frac{\gamma}{1 + \gamma n}  \approx \ds \frac{1}{n} + \frac{1}{\gamma n^2} + O\left( \ds \frac{1}{n^3} \right)
\ee Here is a sketch of Simon's proof: he considered $L_n$, the $(n+1)\times(n+1)$ matrix defined as $(L_n)_{jk}=c_{j-k}$, where $c_{j}=\int e^{- i j \theta} d\mu(\theta)$ is the $j$-th moment of the measure. It is well-known that if $\Phi_n(z)=a_n z^n + a_{n-1} z^{n-1} + \dots + a_{0}$, $\delta_n = (0, 0, \dots, 0, 1)$ and $\langle, \rangle$ being the Euclidean norm,
\be
(a_0, a_1, \dots, a_n) = \left\langle \delta_n, L_n^{-1} \delta_n \right \rangle^{-1} L_n^{-1} \delta_n
\ee

Therefore, the aim is to compute $L_n^{-1}$. By (\ref{dnudef}), $c_n=(1-\gamma) \delta_{n 0} + \gamma$. Let $P_j$ be the $j \times j$ matrix which is $j^{-1}$ times the matrix of all $1$'s, so it is a rank one projection. $L_n$ could be decomposed as
\be
L_n=(1-\gamma) \mathbf{1} + (n+1) \gamma P_{n+1}
\label{ln}
\ee
From (\ref{ln}), one could deduce that the inverse of $L_n$ is
\be
L_n^{-1}=(1-\gamma)^{-1}(\mathbf{1}-P_{n+1})+(1+n\gamma)^{-1} P_{n+1}
\ee

Unfortunately, the method used to prove the result above no longer gives such a nice result when there are two pure points. For instance, we won't have the decomposition as in $(\ref{ln})$, because $L_n$ will be a rank $m$ perturbation of $(1-\sum_{j=1}^{m} \gamma_j) \mathbf{1}$ instead, so the computations will be much more complicated. Besides, this method only works for adding one point to $d\theta/2\pi$ but fails for more general measures. Therefore, we need  another method to attack the problem.

\medskip


From formula (\ref{myformula}) we could make a few observations concerning successive Verblunsky coefficients $\alpha_{n+1}(d\nu)$ and $\alpha_n(d\nu)$: first, we use the fact that $\ol{\vp_{n+1}(\zeta)}=\ol{\zeta^{n+1}} \vp_{n+1}^*(\zeta)$ and rewrite formula (\ref{myformula}) as
\be
\alpha_n(d\nu) = \alpha_n + \ds \frac{(1-|\alpha_n|^2)^{1/2}}{(1-\gamma) \gamma^{-1}+ K_n(\zeta)} \ol{\zeta^{n+1}} \vp_{n+1}^*(\zeta) \vp_n^*(\zeta)
\label{intermediateformula}
\ee

Let $t_n$ be the tail term in the right hand side of (\ref{intermediateformula}) above. Suppose we can prove that $\vp_n^*(\zeta)$ tends to some non-zero limit $L$ as $n$ tends to infinity, then $1/K_n = O(1/n)$, hence,
\begin{equation}
\frac{1} {(1-\gamma)\gamma^{-1} + K_{n}(\zeta)}= \ds \frac{1}{K_n(\zeta)} + O\left( \frac{1}{n^2} \right)
\end{equation}

Besides,  $(\alpha_n)_{n=0}^{\infty}$ is $\ell^2$, therefore $(1-|\alpha_n|^2)^{1/2} \to 1$. As a result,
\be \alpha_n(d\nu) = \alpha_n + t_n \approx \alpha_n + \ds \frac{\ol{\zeta^{n+1}} L^2} {n|L|^2} + \ds o\left( \frac{1}{n} \right)
\ee

Indeed, we shall prove that if $\zeta t_{n+1} - t_n$ is summable, by Theorem \ref{bookthm}, $\lim_{n \to \infty} \vp_n^*(z, d\mu_1)$ exists away from $z=1$. As a result, if we add another a pure point to $d\mu_1$, we can use a similar argument to the one above and formula (\ref{myformula}) to prove that $\alpha_n(d\nu)$ is the sum of $\alpha_n(d\mu_0)$ plus two tail terms and an error term.

In general, if we have a measure $d\mu_m$ as defined in (\ref{dmumdef}), then we add one pure point after the other and use formula (\ref{myformula}) inductively. Therefore, we shall be able to express $\alpha_n (d\mu_m)$ as the sum of $\alpha_n(d\mu_0)$ plus $m$ tail terms, and an error term
\be
\alpha_n (d\mu_m) = \alpha_n (d\mu_0) + t_{1,n} + t_{2,n} + \dots + t_{m,n} + error
\ee By an argument similar to the one above we observe that $t_{j,n}$ is $O(1/n)$ and $z_j t_{j,n} - t_{j,n-1}$ is small. Of course, the 'smallness' has to be determined by rigorous computations that we shall present in the proof Nonetheless, these observations led us to introduce the notion of generalized bounded variation $W_m$, and from that we could deduce that $\lim_{n \to \infty}\vp_n^*(z, d\mu_m)$ exists.

\medskip

\section{Proof of Theorem \ref{theorem0}}\label{proof0}
In the proof of Theorem 10.13.7 in \cite{simon2}, Simon gave the following formula for the Verblunsky coefficients of $d\nu$
\begin{equation}
\alpha_{n}(d\nu)=\alpha_{n}-q_{n}^{-1} \gamma \ol{\vp_{n+1}(\zeta)} \left( \ds \sum_{j=0}^{n} \alpha_{j-1} \frac{\|\Phi_{n+1}\|}{\|\Phi_j\|}\vp_j(\zeta) \right)
\label{mainformula}
\end{equation} where
\begin{eqnarray}
K_{n}(\zeta) & = & \ds \sum_{j=0}^{n} |\vp_j(\zeta)|^2 \label{kndef}\\
q_{n} & = & (1-\gamma) + \gamma K_{n}(\zeta) \\
\alpha_{-1} & = & -1
\end{eqnarray} and all objects without the label $(d\nu)$ are associated with the measure $d\mu$. \\

First, we observe that $\alpha_{j-1}=-\ol{\Phi_j(0)}$, therefore, $\alpha_{j-1}/\|\Phi_j\| = - \ol{\vp_j(0)}$. Second, observe that $\|\Phi_{n+1}\|$ is independent of $j$ so it could be taken out from the summation. As a result, formula (\ref{mainformula}) becomes
\be
\alpha_{n}(d\nu)=\alpha_{n}(d\mu) + q_{n}^{-1} \gamma \ol{\vp_{n+1}(\zeta)} \|\Phi_{n+1}\| \left( \ds \sum_{j=0}^{n} \ol{\vp_j(0)}\vp_j(\zeta) \right)
\label{intermediate}
\ee

Then we use the Christoffel-Darboux formula, which states that for $x, y \in \mathbb{C}$ with $x \ol{y} \not = 1$,
\begin{equation}
(1-\ol{x} y) \left( \ds \sum_{j=0}^{n} \ol{\vp_j(x)} \vp_j(y) \right) = \ol{\vp_{n}^*(x)}\vp_{n}^*(y) -  \ol{x}y \ol{\vp_n(x)} \vp_n (y)
\end{equation}

Besides, note that $q_n^{-1} \gamma = ((1-\gamma)\gamma^{-1} + K_n(\zeta))^{-1}$
As a result, (\ref{intermediate}) could be simplified as follows
\be 
\alpha_n(d\nu)
= \alpha_n + \ds \frac{\ol{\vp_{n+1}(\zeta)} \vp_n^*(0) \vp_n^*(\zeta)}{(1-\gamma)\gamma^{-1} + K_n(\zeta)} \|\Phi_{n+1}\|
\label{formula1}
\ee

Finally, observe that $\vp_n^*(0)=\|\Phi_n\|^{-1}$ and that by (\ref{normphi}), $\|\Phi_{n+1}\|/\|\Phi_n\|=(1-|\alpha_n|^2)^{1/2}$. This completes the proof.

\section{Proof of Theorem \ref{mainthm2}}

The technique used in this proof is a generalization of the one used in proving Theorem \ref{bookthm}. It involves Pr\"ufer variables which are defined as follows

\textbf{Definition} Suppose $z_0 = e^{i\eta} \in \T$ with $\eta \in  [0, 2\pi)$. Define the \textit{Pr\"ufer variables by}
\begin{equation}
\Phi_n(z_0) = R_n(z_0) \exp(i(n \eta + \theta_n(z_0)))
\label{pruefer}
\end{equation} where $\theta_n$ is determined by $|\theta_{n+1} - \theta_n| < \pi$. Here, $R_n(z)=|\Phi_n(z)| >0$, $\theta_n$ is real. By the fact that $\Phi_n^*(z) = z^n \ol{\Phi_n(z)}$ on $\T$, (\ref{pruefer}) is equivalent to
\begin{equation}
\Phi_n^*(z)=R_n(z) \exp(-i\theta_n)
\end{equation}

Under such definition,
\be
\log \left( \ds \frac{\Phi_{n+1}^*}{\Phi_{n}^*} \right) = \log(1 - \alpha_n \exp(i[(n+1)\eta + 2 \theta_n]))
\label{fraclog}
\ee For simplicity, we let
\be a_n = \alpha_n \exp(i[(n+1)\eta + 2 \theta_n])
\label{andef}
\ee

Now write $\log \Phi_{n+1}^*$ as a telescoping sum
\be
\log \Phi_{n+1}^*(z) = \ds \sum_{j=0}^{n} \left( \log \Phi_{j+1}^*(z) - \log \Phi_j^*(z) \right) = \ds \sum_{j=0}^{n} \ds \log \left(\frac{\Phi_{j+1}^*(z)}{\Phi_{j}^*(z)} \right)
\ee
Note that for $|w|\leq Q < 1$, there is a constant $K$ such that
\be
\left| \log(1-w) - w \right| \leq K |w|^2 
\ee Together with (\ref{fraclog}), we have
\be
\log(\Phi_{n+1}^*(z)) = - \ds \sum_{j=0}^{n} \left( a_j + L(a_j) \right)
\label{logphin}
\ee where $|L(a_j)| \leq K |a_j|^2$.

Recall that by assumption, $(\alpha_n(d\mu_0))_{n=0}^{\infty}$ is $\ell^2$. Therefore, by (\ref{andef}), $(a_n)_{n=0}^{\infty}$ is also $\ell^2$, thus $\sum_{j=0}^\infty L(a_j) < \infty$. As a result, in order to prove that $\lim_{n \to \infty} \Phi_{n}^*(z)$ exists, it suffices to prove that $\sum_{j=0}^{\infty} a_j$ exists. \\

Let
\begin{equation}
h_n^{(k)} = \ds \sum_{j=0}^{n-1} \ol{\zeta_k}^j e^{ij\eta} = \ds \frac{\ol{\zeta_k}^n e^{in\eta}-1}{\ol{\zeta_k} e^{i\eta}-1}
\end{equation} Then
\begin{eqnarray}
h_{n+1}^{(k)} - h_n^{(k)} & = & \ol{\zeta_{k}}^n e^{i n \eta} \\
\text{and \quad}|h_n^{(k)} | & \leq & 2| \ol{\zeta_k} e^{i \eta}-1|^{-1}
\end{eqnarray}

Let $g_j = \eta + 2\theta_j$ and recall that $\alpha_n = \sum_{k=1}^{p} \beta_{n,k} $. By rearranging the order of summation, we get
\begin{equation}
S_n = \ds \sum_{j=0}^{n} \alpha_j e^{i(j \eta + g_j)} = \ds \sum_{j=0}^{n} \ds \left( \sum_{k=1}^{p} \beta_{j,k} \right) e^{i(j \eta + g_j)} = \sum_{k=1}^{p} B_n^{(k)}
\label{sndef}
\end{equation} where
\begin{equation}
B_n^{(k)}= \ds \sum_{j=0}^{n} \beta_{j,k} e^{i(j \eta + g_j)}
\end{equation}

We are going to sum by parts by Abel's formula. Suppose $(a_j)_{j=0}^{\infty}$ is a sequence, we define
\begin{eqnarray}
(\delta^+ a)_j & = & a_{j+1} - a_j \\
(\delta^- a)_j & = & a_{j} - a_{j-1}
\end{eqnarray}

Abel's formula states that
\begin{equation}
\ds \sum_{j=0}^{n}(\delta^+ a)_j b_j = a_{n+1} b_n - a_0 b_{-1} - \ds \sum_{j=0}^{n} a_j (\delta^- b)_j
\label{abelformula}
\end{equation}

Now we apply Abel's formula to $B_n^{(k)}$
\begin{equation}
\begin{array}{lll}
B_n^{(k)} & = & \ds \sum_{j=0}^{n}  (\delta^+ h^{(k)})_j  ({\zeta_k}^j \beta_{j,k} e^{i g_j}) \\
& = & h_{n+1}^{(k)} {\zeta_k}^n \beta_{n,k} e^{i g_n} - h_0^{(k)} \zeta_k \beta_{-1,k} e^{i g_{-1}} - \ds \sum_{j=0}^{n} h_j^{(k)} \delta^-( {\zeta_k}^j \beta_{j,k} e^{ig_j})_j
\end{array}
\label{eqn4}
\end{equation}

Note that the term $h_0 \zeta_k^{-1} \beta_{-1,k} e^{i g_{-1}}$ will be canceled in (\ref{eqn4}), without loss of generality we may assume it to be $0$. 


We want to obtain a bound for $B_{n}^{(k)}$. Observe that
\begin{equation}
|\beta_{n,k} | \leq  \sum_{q=1}^{n} |\beta_{q,k} - \ol{\zeta_k} \beta_{q-1,k}| + |\beta_0|  \leq D_k 
\end{equation} where
\be
D_k = \ds \sum_{q=0}^{\infty} |\beta_{q,k} - \ol{\zeta_k} \beta_{q-1,k}| 
\ee is finite because $d\mu \in W_p(\zeta_1, \zeta_2, \dots, \zeta_p)$.

Next, we use the triangle inequality and $|e^{ix}-e^{iy}|\leq |x-y|$ to obtain
\begin{equation}
\begin{array}{lll}
|\delta^-( \zeta_k^j \beta_{j,k} e^{ig_j})_j| & \leq & | \beta_{j,k} (e^{i g_{j}} - e^{i g_{j-1}})| + |{\zeta_k} \beta_{j,k} - \beta_{j-1,k}| \\
& & \\
& \leq & |\beta_{j,k} (\theta_{j} - \theta_{j-1}) | +  |{\zeta_k} \beta_{j,k} - \beta_{j-1,k}|
\end{array}
\end{equation}

It has been proven for Pr\"ufer variables (see Corollary 10.12.2 of \cite{simon2}) that
\begin{equation}
|\theta_{n+1} - \theta_n| < \frac{\pi}{2} |\alpha_n| ( 1 - |\alpha_n|)^{-1}  
\end{equation}

Now recall our assumption that for  $1 \leq k \leq p$, $(\beta_{n,k})_{n=0}^{\infty}$ is $\ell^2$, therefore $\beta_{n,k} \to 0$, $\alpha_n \to 0$, which implies $Q = \sup_n |\alpha_n| < 1$ and $C=\sup_n (1-|\alpha_n|)^{-1} = (1 - Q)^{-1}$. For any $n$ we have
\begin{equation}
|B_{n,k}|  \leq  |\zeta_k e^{i \eta} - 1|^{-1} \left( 2 D_k + \ds \frac{\pi}{2}  \sum_{j=0}^{\infty} |\beta_{j+1,k}| |\alpha_j| (1 - Q)^{-1} \right) < \infty
\end{equation} It follows that $\sup_{n} |S_n|<\infty$. This proves (\ref{mainthm21}). \\

The computations above also show that the sum in the right hand side of (\ref{eqn4}) is absolutely convergent as $n \to \infty$ and the convergence is uniform on any compact subset of $\T \backslash \{ \zeta_1, \zeta_2, \dots, \zeta_p\}$. Therefore, $\lim_{j \to \infty} \beta_{j,k}=0$ for all $1 \leq k \leq p$ implies that $ \lim_{n \to \infty} B_{n,k}$ exists, thus $\lim_{n \to \infty} S_n$ exists and is finite. This proves (\ref{mainthm22}). \\

Moreover, for each fixed $k$, $(\beta_{n,k})_{n=0}^{\infty}$ is $\ell^2$, $(\alpha_n)_{n=0}^{\infty}$ is also $\ell^2$, hence the Szeg\"o function $D(z)$ exists and it has boundary values a.e.. Now decompose $d\mu = w(\theta)\frac{d\theta}{2\pi} + d\mu_s$. It is well-known that $\Phi_n^* \to D(0) D^{-1}$ in $L^2(w(\theta)\frac{d\theta}{2\pi})$. Since $\Phi_n^* \to \Phi_\infty^*$ uniformly on any compact subset of $\T \backslash \{ \zeta_1, \zeta_2, \dots, \zeta_p\}$, the limit also converges in the $L^2$-sense. Besides, it is well known that $D(0)=\lim_{n \to \infty}\|\Phi_n\|=\prod_{n=0}^{\infty}(1-|\alpha_n|^2)^{1/2}$, hence
\begin{eqnarray}
\Phi_\infty^*(z) & = & D(0) D^{-1} (z) \label{lhs} \\
\vp_\infty^*(z) & = & D^{-1} (z)
\end{eqnarray} on $\T \backslash \{\zeta_1, \zeta_2, \dots, \zeta_p \}$.  \\

\medskip

\section{Proof of Theorem \ref{mainthm1}}




We proceed by induction.

\subsection{Base Case} \label{basecase} Let any object without the label $(d\mu_1)$ be associated with the measure $d\mu_0$. First we start by considering adding one pure point $z_1 = e^{i \omega_1} \in \T$, $\omega_1 \not = 1$, to $d\mu_0 \in W_1(1)$ which has $\ell^2$ Verblunsky coefficients.

Define $\tilde{\xi}_n(d\mu_1)$ as
\begin{equation}
\tilde{\xi}_n (d\mu_1)= \ds \frac{(1-|\alpha_n|^2)^{1/2}}{(1-\gamma) \gamma^{-1}+ K_n(z_1)} \ol{\vp_{n+1}(z_1)} \vp_n^*(z_1)
\label{equation5}
\end{equation} where $\alpha_j = \alpha_j (d\mu_0)$ and $(\Phi_n)_{n=0}^{\infty}$ is the family of orthogonal polynomials for $d\mu_0$. Because of formula (\ref{myformula}), we want to simplify $\tilde{\xi}_n(d\mu_0)$.

Since $d\mu_0 \in W_{1}(1)$ and $\sum_{j=0}^{\infty} |\alpha_j|^2 < \infty$, by Theorem \ref{mainthm2} $\lim_{n \to \infty} \vp_n^*(z_1)= D(z_1)^{-1}$, which implies $1/K_n(z_1) = O(1/n)$. Hence,
\be
\tilde{\xi}_n (d\mu_1)= \ds \frac{(1-|\alpha_n|^2)^{1/2}}{K_n(z_1)} \ol{\vp_{n+1}(z_1)} \vp_n^*(z_1) + \ds \bigonsq
\label{36}
\ee

Moreover, $\ol{\vp_{n+1}(z_1)} = \ol{z_1^{n+1}} \vp_{n+1}^*(z_1)$. We can further simplify and obtain
\be
\ds \alpha_n(d\mu_1) = \alpha_n + \ol{z_1}^{n+1} \ds 
\frac{D(z_1)^{-2}}{|D(z_1)|^{-2}} \frac{1}{n}+ o\left( \ds \frac{1}{n}\right)
\label{form1}
\ee

Let $c_1 = \ol{z_1} D(z_1)^{2}/ |D(z_1)|^2$. This proves (\ref{equation1}) for $m=1$.

\medskip

\textbf{Remark:} Note that the error term in the right hand side of (\ref{form1}) is dependent on $\gamma_1$. This is because as $\gamma_0 \to 0$, $d\mu_1 \rightarrow d\mu_0$ weakly, which implies that for each $n$, $\alpha_n(d\mu_1) \to \alpha_n(d\mu_0)$. Since the tail term $\ol{z_1}^{n+1} \frac{D(z_1)^{-2}}{|D(z_1)|^{-2}} \frac{1}{n}$ in (\ref{form1}) is independent of $\gamma_1$, if the error term is also independent of $\gamma_1$, $\alpha_n(d\mu_1) \not \to \alpha_n(d\mu_0)$.

\medskip
It remains to show the claimed properties of $\Phi_n(d\mu_1)$. To do that, it suffices to show that $(\alpha_n(d\mu_1))_{n=0}^{\infty}$ is $\ell^2$ and it is in the class $W_{2}(1, z_1)$, then we can conclude by Theorem \ref{mainthm2}.

First of all, it is clear that $(\alpha_n(d\mu_1))_{n=0}^{\infty}$ is $\ell^2$ because $(\alpha_n)_{n=0}^{\infty}$ is $\ell^2$ and $\tilde{\xi_n}(d\mu_1)$ is $O(1/n)$.

Next, we want to show that
\be
\ds \sum_{n=0}^{\infty} |z_1 \tilde{\xi}_{n+1} - \tilde{\xi}_n | < \infty
\label{difference}
\ee By (\ref{36}), the error term is in the order of $O(1/n^2)$, therefore this is the same as showing the following is $\ell^1$-summable
\begin{equation}
\left| \ds \frac{\vp_{n+2}^*(z_1)\vp_{n+1}^*(z_1) (1-|\alpha_{n+1}|^2)^{1/2}}{K_{n+1}} - \ds \frac{\vp_{n+1}^*(z_1)\vp_{n}^*(z_1) (1-|\alpha_{n}|^2)^{1/2}}{K_{n}} \right|
\end{equation}

We are going to estimate term by term.

\begin{itemize}
\item Let $\rho_n = (1-|\alpha_n|^2)^{1/2}$. We estimate the following using the recurrence relation for orthogonal polynomials (\ref{normrec})
\begin{equation}
\begin{array}{ll}
& \vp_{n+1}^*(z_1) - \vp_{n}^*(z_1) \\
= & (\rho_n \vp_n^*(z_1) - \alpha_n \vp_{n+1}(z_1)) - \vp_n^*(z_1) \\
= & (\rho_n -1) \vp_n^*(z_1) - \alpha_n \vp_{n+1}(z_1)
\end{array}
\end{equation}

Since $\rho_n - 1 = O(|\alpha_n|^2)$, $\vp_n^*(z_1) = D(z_1)^{-1} + o(1)$ and $1/K_n = O(1/n)$,
\begin{equation}
\left| \vp_{n+1}^*(z_1) - \vp_{n}^*(z_1) \right| = (O(|\alpha_n|^2) + |\alpha_n|)|D(z_1)^{-1} + o(1)| = O(|\alpha_n|)
\end{equation} Hence,
\begin{equation}
\left| \ds \frac{\left(\vp_{n+1}^*(z_1) - \vp_{n}^*(z_1) \right)\vp_{n+1}^*(z_1) (1-|\alpha_n|^2)^{1/2}}{K_n} \right| = O \left( \ds \frac{|\alpha_n|}{n}\right)
\label{esteqn1}
\end{equation}

\item
If we change $n$ to $n+1$, the same argument still holds. Therefore, 
\begin{equation}
\left| \ds \frac{\left(\vp_{n+2}^*(z_1) - \vp_{n+1}^*(z_1) \right)\vp_{n}^*(z_1) (1-|\alpha_{n}|^2)^{1/2}}{K_{n}} \right| = O \left( \ds \frac{|\alpha_{n+1}|}{n}\right)
\label{esteqn2}
\end{equation}

\item Observe that
\begin{equation}
|(1-|\alpha_{n+1}|)^{1/2} - (1-|\alpha_{n}|)^{1/2}|= O(|\alpha_n| + |\alpha_{n+1}|)
\end{equation} Hence,
\begin{equation}
\ds \left| \frac{\left[(1-|\alpha_{n+1}|)^{1/2} - (1-|\alpha_{n}|)^{1/2} \right] \vp_{n+1}^*(z_1)\vp_{n}^*(z_1)}{K_n} \right| = O\left( \ds \frac{|\alpha_{n+1}| + |\alpha_{n}|}{n} \right)
\label{esteqn3}
\end{equation}

\item Finally, note that
\begin{equation}
\left( \ds \frac{1}{K_{n+1}} - \frac{1}{K_{n}} \right) \vp_{n+1}^*(z_1)\vp_{n}^*(z_1)(1-|\alpha_n|^2)^{1/2} = \bigonsq 
\label{esteqn4}
\end{equation}

\end{itemize}
\medskip

Combining all the estimates above, we have
\begin{equation}
|z_1 \tilde{\xi}_{n+1} - \tilde{\xi}_n | = O\left( \ds \frac{|\alpha_n|+|\alpha_{n+1}|}{n} \right) + \bigonsq
\label{tildexidiff}
\end{equation}

As a result,
\begin{equation}
\ds \sum_{n=0}^{\infty} |z_1 \tilde{\xi}_{n+1} - \tilde{\xi}_n | < \infty
\end{equation} and by Theorem \ref{mainthm2}, the proof of the case $m=1$ is complete.

\medskip

\subsection{Induction Step} \label{induction step} We consider $d\mu_m$ as defined in (\ref{dmumdef}) as a measure formed by adding a pure point to $d\mu_{m-1}$ in the following manner

Let
\begin{equation}
\tilde{\gamma_j} = (1- \gamma_m)^{-1} \gamma_j
\end{equation} and
\begin{equation}
d\mu_{m-1} =\left(1 - \ds \sum_{l=1}^{m-1} \tilde{\gamma_l} \right) d\mu_0 + \ds \sum_{l=0}^{m-1} \tilde{\gamma_l} \delta_{\omega_l}
\label{dmum1def}
\end{equation}
Then we could write
\begin{equation}
d\mu_m = (1- \gamma_m) d\mu_{m-1}+ \gamma_m \delta_{\omega_m}
\end{equation}

Recall that $0<\sum_{l=1}^{m} \gamma_l <1$, or equivalently, $\sum_{l=1}^{m-1} \gamma_l < 1- \gamma_m$. Hence, 
\begin{equation}
0< \ds \sum_{j=1}^{m-1} \tilde{\gamma_j} = \left( 1- \gamma_m \right)^{-1}  \left( \ds \sum_{j=1}^{m-1} \gamma_j \right) < 1
\end{equation}

Therefore, $d\mu_{m-1}$ satisfies the induction hypothesis, so its family of Verblunsky coefficients is $\ell^2$ and $d\mu_{m-1} \in W_{m}(1, z_1, z_2, \dots, z_{m-1})$. Hence, $\lim_{n \to \infty} \vp_n^*(z_m, d\mu_{m-1})$ exists and is equal to $(1-\sum_{j=1}^{m-1}\gamma_j)^{1/2} D(z_m, d\mu_0)^{-1}$ (see remark following Theorem \ref{mainthm1}). As a result, we can use a similar argument as in the base case and deduce that
\be \begin{array}{ll}
\alpha_n(d\mu_m) & = \alpha_n(d\mu_{m-1}) + \ds \ol{z_m}^{n+1} \frac{|D(z_m, d\mu_0)|^{2}}{D(z_m, d\mu_0)^{2}} \frac 1 n + E_n \\
\\
& = \alpha_n (d\mu_0) + \ds \sum_{j=1}^{m} \ds \frac{\ol{z_j}^{n} c_j}{n} + E_n
\end{array}
\ee where $c_j = \ol{z_j} D(z_j, d\mu_0)^{2}/|D(z_j, d\mu_0)|^2$, $1 \leq j \leq m$, are constants independent of the weights $\gamma_1, \gamma_2, \dots, \gamma_m$ and of $n$; and $E_n = E_n(z_1, z_2, \dots, z_m, \gamma_1, \gamma_2, \dots, \gamma_m)$ is in the order of $o(1/n)$. This proves (\ref{equation1}).

By estimating consecutive Verblunsky coefficients in the same way we did in the base case, we prove that $d\mu_m \in W_{m+1} (1, z_1, z_2, \dots, z_m)$. Thus, we can apply Theorem \ref{mainthm2} to prove that $\vp_n^*(z_m)$ tends to $D(z_m, d\mu_m)^{-1}$. This completes the proof of Theorem \ref{mainthm1}. \medskip

\textbf{Remark:} Note that if $d\mu_0 \in W_p(\zeta_1, \zeta_2, \dots, \zeta_p)$ and $z_j \not = \zeta_k$ for all $j,k$, we can use the same arguments as in the proof of Theorem \ref{mainthm1} to prove similar results, i.e., $\alpha_n(d\mu_m)$ is in the form (\ref{equation1}), $d\mu_m$ is in $W_{m+p}(\zeta_1, \zeta_2, \dots, \zeta_p, z_1, z_2, \dots, z_m)$ and that $\lim_{n \to \infty} \vp_{n}(z,d\mu_m)=D(z, d\mu_m)^{-1}$ for $z \not = \zeta_1, \zeta_2, \dots, \zeta_p, z_1, z_2, \dots, z_m$.

\bigskip
\section{Acknowledgements}
I would like to thank Professor Barry Simon for suggesting this problem and proof-reading this paper, as well as his patience and enthusiasm for advising his students; last but not least, for writing the two great reference books \cite{simon1, simon2}. \\

\end{document}